\documentclass[
final
]{dmtcs-episciences}

\usepackage[utf8]{inputenc}
\usepackage{subfigure}
\usepackage{amssymb}
\usepackage{amsthm}
\usepackage{graphics}
\usepackage{epsfig}
\usepackage{amssymb}
\usepackage{graphicx}
\usepackage{subfigure}
\usepackage{color}
\usepackage{tikz}
\usepackage{tikz,pgfplots}
\pgfplotsset{compat=1.12,axis lines=center}
\usepgfplotslibrary{fillbetween}
\usetikzlibrary{patterns}
\usepackage{hyperref}
\usepackage{a4wide}
\usepackage{amsmath}
\usepackage{amsfonts}
\usepackage{enumerate}
\usepackage{multicol}

\newcommand{\pf}{\noindent {\bf Proof: }}
\usepackage[round]{natbib}

\newtheorem{lemma}{\bf Lemma}[section]

\newtheorem{theorem}{\bf Theorem}[section]
\newtheorem{proposition}[lemma]{\bf Proposition}
\newtheorem{corollary}[lemma]{\bf Corollary}
\newtheorem{definition}{\bf Definition}[section]
\newtheorem{remark}{\bf Remark}[section]

\author{Angsuman Das\affiliationmark{1}\thanks{Supported by DST grants Sanction no. SRG/2019/ 000475 and SR/FST/MS-I/2019/41, Govt. of India.}
  \and Hiranya Kishore Dey\affiliationmark{2}\thanks{Supported by CSIR-SPM Fellowship, IIT Bombay Post Doctoral fellowship and Post Doctoral Fellowship, HRI Prayagraj, Dept. of Atomic Energy, Govt. of India.}}
\title[Determining Number of Kneser Graphs]{Determining Number of Kneser Graphs: Exact Values and Improved Bounds}
\affiliation{
  Department of Mathematics, Presidency University, Kolkata, India\\
  Department of Mathematics, Harish-Chandra Research Institute, Prayagraj (Allahabad) India}
\received{2021-06-29}
\revised{2021-11-11}
\accepted{2022-02-16}
\begin{document}
\publicationdetails{24}{2022}{1}{10}{7627}
\maketitle
\begin{abstract}
  The determining number of a graph $G = (V,E)$ is the minimum cardinality of a set $S\subseteq V$ such that pointwise stabilizer of $S$ under the action of $Aut(G)$ is trivial. In this paper, we provide some improved upper and lower bounds on the determining number of Kneser graphs. Moreover, we provide the exact value of the determining number for some subfamilies of Kneser graphs. 
\end{abstract}

\section{Introduction}

The {\it determining number}, denoted by $Det(G)$, of a graph $G = (V,E)$ is the minimum cardinality of a set $S\subseteq V$ such that the automorphism group of the graph obtained from $G$ by fixing every vertex in $S$ is trivial and such a set $S$ is called a determining set. It was introduced independently by \cite{boutin} and (defined as {\it fixing number}) \cite{harary} in 2006 as a measure of destroying the symmetry of a graph. Determining sets are quite useful in investigating graph automorphisms and are even more useful in proving that all automorphisms of a given graph have been found. The authors in \cite[Proposition 2]{boutin} linked the size of a determining set to the size of the automorphism group and proved that the determining number of a graph is greater than equal to the logarithm of the cardinality of the automorphism group. 
Apart from that, determining sets have also been used to find distinguishing number. For reference, one can see (\cite{boutin-distinguish}). In that paper the authors used determining sets to show that the Kneser graphs $K(n,k)$ with $ n \geq 6$ and $k \geq 2$ are 2-distinguishable graphs. Apart from proving general bounds and other results on determining number, researchers have attempted to find exact values of determining number of various families of graphs like Kneser Graphs (\cite{kneser}), Coprime graphs (\cite{coprime}), Generalized Petersen graphs (\cite{das-caldam}), \cite{das-saha} etc. 

\subsection{Preliminaries}
The agenda of finding the determining sets of Kneser graphs was initiated in the introductory paper by \cite{boutin}. The next attempt towards it was done in \cite{kneser}.

The Kneser graph $K(n,k)$ has vertices associated with the $k$-subsets of the $n$-set $[n] = \{1, \ldots , n\}$ and edges connecting disjoint sets. This family of graphs is usually considered for $n \geq 2k$ but here we shall assume that $n > 2k$ since the case $n = 2k$ gives a set of disconnected edges and its determining number is half the number of vertices. It is known that the automorphism group of $K(n,k)$ is isomorphic to $S_n$. The following results regarding determining sets of Kneser graphs were proved in \cite{boutin} and \cite{kneser}.

{\lemma \label{vital} \cite{boutin} The set $S=\{V_1 , \ldots , V_r \}$ is a determining set for $K(n,k)$ if and only if there exists no pair of distinct elements $a, b \in [n]$ so that for each $i$ either $\{a, b\} \subseteq V_i$ or $\{a, b\} \subseteq V^c_i$.}

{\remark Thus $S$ is a determining set if for all $a,b \in [n]$, there exists $i$ such that $a \in V_i$ and $b\not\in V_i$. In that case, we say that $V_i$ separates $a$ and $b$. In other words, it tells that
	$V_1,V_2,\ldots,V_r$ separates  any pair $a,b \in [n]$ and 
	$\displaystyle\cup_{i=1}^{r} V_i$ can miss at most one element of $[n]$. }

{\proposition \label{weak-bound} \cite{boutin},\cite{kneser} $\left\lceil\log_2(n+1)\right\rceil\leq Det(K(n,k))\leq n-k$. }

\cite{boutin} also showed that $Det(K(2^r-1,2^{r-1}-1))=r$. In \cite{kneser}, the authors associate a $k$-regular hypergraph with every subset of vertices of $K(n,k)$. The determining set is achieved by imposing conditions on the edges of the corresponding hypergraph. Using this idea of hypergraphs, they proved the following theorems. 

{\theorem \label{caceres_first_main_result} Let $k$ and $d$ be two positive integers such that $k \leq d$ and $d > 2$. Then $$Det\left(K\left(\left\lfloor \dfrac{d(k+1)}{2}\right\rfloor +1,k \right)\right)=d.$$ } 

{\theorem \label{thm:caceres_second_main_result} Let $k$ and $d$ be two positive integers where $3 \leq k + 1 \leq d$. For every $n\in \mathbb{N}$ such that $\left\lfloor \dfrac{(d-1)(k+1)}{2}\right\rfloor < n  < \left\lfloor \dfrac{d(k+1)}{2}\right\rfloor$, it holds that $Det(K(n+1,k))=d$. }

Caceres {\it et.al.} also answered the following question, posed in \cite[Question 2]{boutin}. 

{\theorem \label{caceres_third_main_result} $Det(K(n,k))=n-k$ if and only if either $k=1$ or $k=2$ and $n=4,5$.\qed}\\

%

\subsection{Our Contribution}
\begin{figure}[ht]
	\centering
	\begin{center}
		\begin{tikzpicture}[scale=1.15]
		
		\draw[->, ultra thick] (0,0)--(12,0) node[right]{};
		\draw[->, ultra thick] (0,0)--(0,9) node[above]{};
		
		\draw[-,blue, line width=2pt] (0.1,0)--(7.05,8);
		\draw[-,black, line width=2pt] (0.8,0)--(7,7.1);
		
		\begin{axis}[samples=100]
		
		\addplot[scale=1,ultra thick,color=blue,name path=A,domain=0:55] plot ({1.4*\x+1},{\x});
		
		\addplot[scale=1,name path=X,domain=0:106] plot ({1*\x},{0});
		\addplot[scale=1,name path=B,domain=0:49]  plot ({1.4*\x+12},{\x});
		
		\addplot[scale=1,name path=C,domain=1.2:14] plot ({\x*(\x+1)/2+1},{\x});
		\addplot[scale=1,name path=D,domain=11:20] plot ({(\x-6)*(\x-5)/2+1},{\x});
		\addplot[scale=1,name path=E,domain=21.1:28] plot ({(\x-14)*(\x-13)/2+1},{\x});
		\addplot[scale=1,name path=F,domain=29.7:35] plot ({(\x-21)*(\x-20)/2+1},{\x});
		\addplot[scale=1,name path=G,domain=37.8:42] plot ({(\x-28)*(\x-27)/2+1},{\x});
		
		\addplot [thick,color=blue,	fill=yellow!100, 	fill opacity=0.4	]
		fill between[of=C and X,soft clip={domain=0:160}	];
		\addplot [thick,color=blue,	fill=gray, 	fill opacity=0.8	]
		fill between[of=C and D	];
		\addplot [thick,color=blue,	fill=gray, 	fill opacity=0.5	]
		fill between[of=E and D	];
		\addplot [thick,color=blue,	fill=gray, 	fill opacity=0.3	]
		fill between[of=E and F	];
		\addplot [thick,color=blue,	fill=gray, 	fill opacity=0.1	]
		fill between[of=G and F	];
		\end{axis}

		\node[rotate=50] at (3.5,4.5) {$n=2k+1 $};
		\node[rotate=48] at (5.4,4.8) {$n=2k+2$};
		\node at (5,.5) { $n\geq \frac{k(k+1)}{2}+1$};
		\node at (9.6,.5) { $\rightarrow Det(K(n,k))$ is known};
		
		\node[rotate=8] at (4,1.4) { $\frac{k(k-1)}{2}+1\leq n\leq \frac{k(k+1)}{2}$};
		\node at (9.2,1.7) { $\rightarrow Det(K(n,k))\leq k$};
		
		\node[rotate=8] at (4.8,2.2) { $\frac{(k-1)(k-2)}{2}+1\leq n\leq \frac{k(k-1)}{2}$};
		\node at (9.5,2.5) { $\rightarrow Det(K(n,k))\leq k-1$};
		
		\node[rotate=8] at (4.8,3.1) {\tiny $\frac{(k-2)(k-3)}{2}+1\leq n\leq \frac{(k-1)(k-2)}{2}$};
		\node at (9.5,3.3) { $\rightarrow Det(K(n,k))\leq k-2$};
		\node at (9.1,4.4) {\Large $\vdots$};
		\node at (9.1,5) {\Large $\vdots$};
		\node at (9.5,8) { $~~~~~~~~~\rightarrow Det(K(n,k))=\lceil \log_2(n+1) \rceil$};
		\node at (9.5,7.15) { $~~~~~~~~~\rightarrow Det(K(n,k))=\lceil \log_2(n+1) \rceil$};
		\node at (9.5,6.5) { ~~~~~~~ ~~~~~or $\lceil \log_2(n+1) \rceil+1$};
		\node[rotate=48] at (3.5,4.5) {};
		\fill[blue] (1,1.1) rectangle (1.2,1.3);
		\fill[blue] (2.5,2.8) rectangle (2.7,3);
		\fill[blue] (4.7,5.3) rectangle (4.9,5.5);
		
		\fill[black] (1,.3) rectangle (1.2,.5);
		\fill[black] (2.5,2) rectangle (2.7,2.2);
		\fill[black] (4.7,4.5) rectangle (4.9,4.7);
		
		\fill[white] (-.15,.5) rectangle (-.65,5.5);
		\fill[white] (1,-.15) rectangle (7,-.5);
		
		\node at (-.6,3.5) {\Large $k$};
		\node at (-.6,4.5) {\Huge $\uparrow$};
		\node at (6.5,-0.5) {\Large $n$};
		\node at (7.5,-0.5) {\Huge $\rightarrow$};
		\end{tikzpicture}
		\caption{Diagramatical representation of exact value and upper bounds of $Det(K(n,k))$}
		\label{diagram}
	\end{center}
\end{figure}
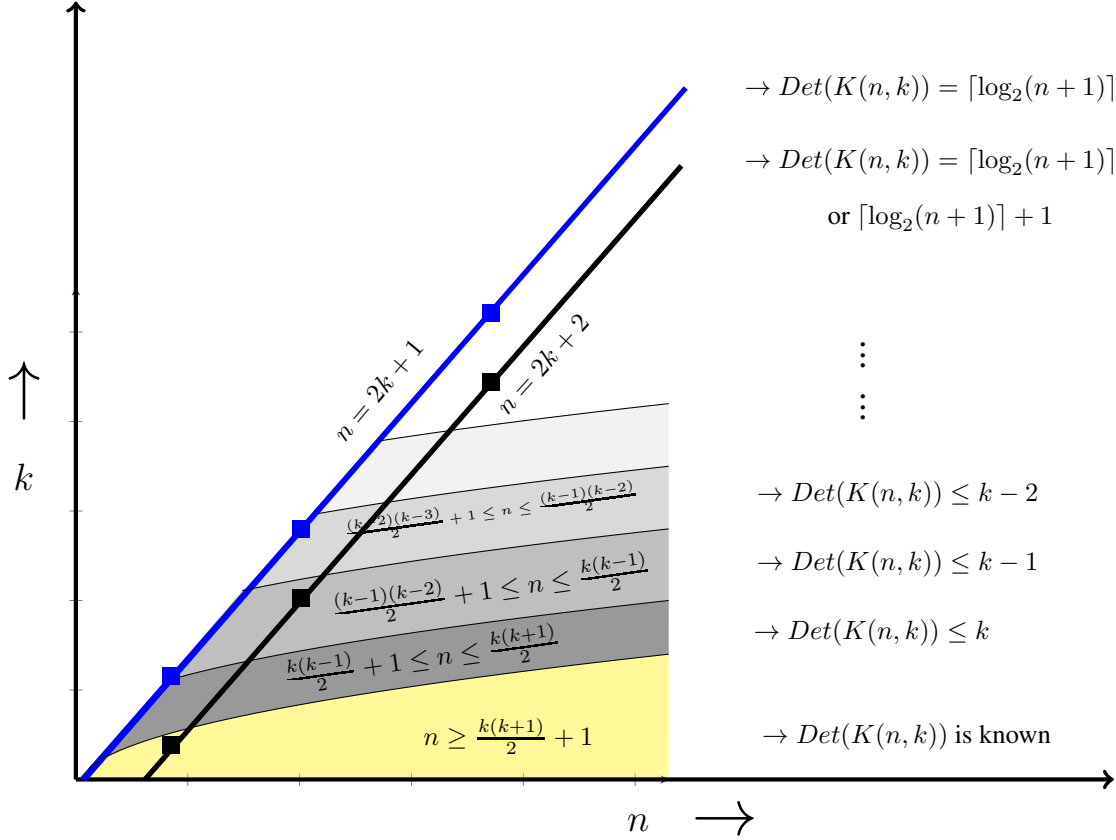

In this paper, we provide some improved bounds on the determining number of Kneser graphs. Moreover, we provide the exact value of the determining number for some subfamilies of Kneser graphs. In Section \ref{recursion-section}, we prove some recursions involving $Det(K(n,k))$ with respect to both $n$ and $k$. In particular, Theorem \ref{thm:rec_for_k_fixed} and Theorem \ref{thm:rec_for_n_fixed} are crucial for the proofs of the main results in the forthcoming sections. In Section \ref{exact-value-section}, we find the exact value of determining number of $K(n,k)$ when $n=2k+1$. It is to be noted that exact value of determining number of $K(n,k)$ for $n=2k+1$ was known earlier only when $n+1$ is a power of $2$ (shown by blue squares in Figure \ref{diagram}, whereas we found the exact value for all the points on the blue straight line in Figure \ref{diagram} (See Theorem \ref{thm:main_theorem_sec2})). We also find $Det(K(n,k))$ when $n=2k+2$ is a power of $2$ (shown by black squares in Figure \ref{diagram}). In fact, we prove that the value of $Det(K(n,k))$ when $n=2k+2$, i.e., $(n,k)$ lie on the bold black straight line shown in Figure \ref{diagram} (even if it is not a power of $2$), is one of the two consecutive integers. In Section \ref{Bounds-section}, we prove some improved lower and upper bounds on $Det(K(n,k))$ compared to those in (\cite{boutin}) and (\cite{kneser}). In particular, Theorem \ref{thn:lowerbound} improves upon the lower bound given in \cite{boutin} and Theorem \ref{thm:stronger_upper_bound} is an improvement of Theorem \ref{thm:Caceres} (Theorem 3.1, \cite{kneser}). Figure \ref{diagram} illustrates the value of $n$ and $k$ for which the exact value of $Det(K(n,k))$ is known or an upper bound is known. The yellow region is where the exact value of determining number was previously known from Theorem \ref{caceres_first_main_result}. Caceres {\it et.al.} also showed that above the yellow region, $Det(K(n,k))\leq k$. We prove a stronger upper bound as shown in different shades of grey in Figure \ref{diagram}, i.e., as we move from darker shades of grey to lighter shades of grey, our upper bounds are progressively tighter than that proved in \cite{kneser}. 
For definitions and terms used in the paper, readers are referred to the classic book \cite{godsil-royle}. 

\section{Recursions}\label{recursion-section}

In this section, we first prove some recursions which we will use throughout the paper. We start with the following result which shows that if we keep $k$ fixed and increase $n$ by $1$, then the determining number either remains same or can at most increase by $1$.  

{\theorem
	\label{thm:rec_for_k_fixed} For all positive integers $n,k$ with $2k<n$, $$Det(K(n,k))\leq Det(K(n+1,k)) \leq Det(K(n,k))+1.$$}
\pf We first prove the first inequality. Let $Det(K(n+1,k))$ be $r$ and let $\{A_1,A_2,\ldots,A_r\}$ be a determining set for $K(n+1,k)$. If $\cup_{i=1}^r A_i$ misses any element in $\{1,2,\ldots,n+1\}$, then by suitable relabelling, $\{A_1,A_2,\ldots,A_r\}$ is a determining set for the graph $(K(n,k))$. Thus, we assume that $\cup_{i=1}^r A_i=[n+1]$. Let $N_i=|\{A_j:i \in A_j\}|$ for $i=1,2\ldots,n+1$. We observe that if $N_i \neq N_j$, then the elements $i$ and $j$ are of course separated. By suitable relabelling, without loss of generality, we can assume that $N_1\geq N_2\geq \cdots \geq N_{n+1}>0$. The idea is to replace $n+1$ in $A_i$'s by some other elements such that $\{A_1,A_2,\ldots,A_r\}$ remains a determining set for $K(n+1,k)$. We start with $A_1$.

If $n+1 \not\in A_1$, we do not change $A_1$. So, let $n+1 \in A_1$. If $1 \not\in A_1$, then replace $n+1$ by $1$ in $A_1$, i.e., $A_1'=A_1\cup \{1\}\setminus \{n+1\}$. Note that this manipulation, increases $N_1$ by $1$, decreases $N_{n+1}$ by $1$ and keeping all other $N_i$'s unchanged, i.e., $N_1'>N_2\geq N_3 \geq \cdots \geq N_n > N_{n+1}'$.  Thus for any $2 \leq i \leq n$,  $1$ and $i$ are of course separated as $N_1' > N_i$ and
$i$ and $n+1$ is also separated as $N_{i} > N_{n+1}$. Moreover, this manipulation does not effect any other pairs $i$ and $j$ with both $2 \leq i \neq j \leq n $ and therefore they are also separated by $\{A_1', A_2, \dots, A_r\}$. Therefore, $\{A_1',A_2,\ldots,A_r \} $ is a determining set for $K(n+1,k)$. Hence we can assume $1 \in A_1$. 

{\bf Claim 1:} $t \notin A_1 \implies $ either $n+1$ can be replaced by $t$ in $A_1$ and
$\{A_1 \setminus  \{n+1\}\cup\{t\}, A_2,\ldots,A_r\}$ still form a determining set of $K(n+1,k)$, or there exists $j <t$ such that $ j \in A_1$ and the couple $(j,t)$ is separated only by $A_1$.

{\bf Proof of Claim 1:} Let $t \notin A_1$. If $n+1$ can be replaced by $t$ and $\{A_1 \setminus  \{n+1\}\cup\{t\}, A_2,\ldots,A_r\}$ is still a determining set of $K(n+1,k)$, then the claim is true. If $t$ can not replace $n+1$, then there must exist some $j\in A_1$ such that $(j,t)$ is separated only by $A_1$. Now, if $t<j$, then $N_t\geq N_j$. Thus, if $t$ replaces $n+1$ in $A_1$, then $N_t'=N_t+1>N_j$. That means there is more sets in $\{A_1',A_2,\ldots,A_r\}$ which contains $t$ than which contains $j$. Thus $t$ and $j$ can be separated by $\{A_1',A_2,\ldots,A_r\}$, a contradiction. Hence $j<t$ and the claim follows.

Hence, for every $t \notin A_1$ such that replacing $n+1$ by $t$ would create a {\it problem}, there exists a $j\in A_1$ with $j<t$. By {\it problem}, we mean a situation where there exist two elements $j$ and $t$ which are not separated. 

{\bf Claim 2:} If $t_1,t_2 \not\in A_1$ with $t_1\neq t_2$ be such that replacing $n+1$ by $t_1$ or $t_2$ is a {\it problem}, then there exists $j_1,j_2 \in A_1$ with $j_1<t_1$, $j_2<t_2$ and $j_1\neq j_2$.

{\bf Proof of Claim 2:} The existence of such $j_1$ and $j_2$ are guaranteed by Claim 1. Only thing left to be shown is that $j_1\neq j_2$. If possible, let $j_1=j_2=j\in A_1$ (say). Thus the pairs $(j,t_1)$ and $(j,t_2)$ are separated only by $A_1$. Thus $(j,t_1)$ (and similarly $(j,t_2)$) are either both present or both absent in $A_2,A_3,\ldots,A_r$. Thus $j,t_1,t_2$ are either all present or all absent in the sets $A_2,A_3,\ldots,A_r$. In particular, $t_1,t_2$ 
\begin{itemize}
	\item are both present or both absent in $A_2,A_3,\ldots,A_r$, and
	\item are both absent in $A_1$.
\end{itemize}
This contradicts the fact that $A_1,A_2,\ldots,A_r$ separates $t_1$ and $t_2$. Hence $j_1\neq j_2$ and the claim follows.

Now as $|A_1|=k$, there exists $n-k+1$ elements in $[n+1]$ which are not in $A_1$. From Claim 1 and 2, either $n+1$ can be replaced by some element in $[n+1]\setminus A_1$ or we get $n-k+1$ distinct elements $j_1,j_2,\ldots,j_{n-k+1}$ in $A_1$. However as $n-k+1>k$, this is a contradiction. Thus $n+1$ can be replaced by some element $t\in [n+1]\setminus A_1$ and $\{A_1'=A_1 \setminus  \{n+1\}\cup\{t\}, A_2,\ldots,A_r\}$ is still a determining set of $K(n+1,k)$.

Thus it is possible to replace $n+1$ in $A_1$. It is to be noted that once $n+1$ is replaced by some $t$ in $A_1$, $N_t$ is increased by $1$, i.e., $N_t'=N_t+1$ and $N_{n+1}$ is decreased by $1$, i.e., $N_{n+1}'=N_{n+1}-1$. Thus, after $A_1$ is modified, in the new sequence of $\{N_i\}$, $N_{n+1}'$ remains the least element. (Note that the ordering of the sequence may change for the term $N_t'$.) We rearrange the terms in the new sequence $\{N_i'\}$ in descending order where $N_{n+1}'=N_{n+1}-1$ remains the smallest term. Also note that, as we relabel the elements in $[n+1]$ to get $N_1'\geq N_2'\geq \cdots \geq N_n'> N_{n+1}'=N_{n+1}-1$, the element $n+1$ is not relabelled as already $N_{n+1}'=N_{n+1}-1$ is the least among the $N_i'$'s. Thus, the element $n+1$ can not re-enter the modified $A_1$ by relabelling. Now, we apply the same process on $A_2$ to get rid of $n+1$ and so on. Continuing in this manner, we replace $n+1$ from each of $A_1,A_2,\ldots,A_r$. Thus $\{A_1,A_2,\ldots,A_r\}$ is a determining set of $K(n,k)$ and hence the theorem.

We now prove the second inequality. Let $V=\{V_1,V_2,\ldots,V_r\}$ be a minimal determining set for $K(n,k)$. Then, $\cup _{i=1}^r V_i$ can miss at most one element of $[n]$. If $\cup _{i=1}^r V_i=[n]$, then $V$ is also a determining set for $K(n+1,k)$. If $\cup _{i=1}^r V_i$ misses one element of $[n]$, without loss of generality, we assume that element to be $n$, i.e., $\cup _{i=1}^r V_i=[n-1]$. We take $V_{r+1}=\{1,2\ldots,k-1,n+1\}$. 
To check that $V'=\{V_1,V_2,\ldots,V_r,V_{r+1}\}$ is a determining set for $K(n+1,k) $, we consider the following cases.
\begin{enumerate}
	\item Let $1 \leq i,j \leq n$. Then $i$ and $j$ are separated by some $V_t$ for $1 \leq t \leq r$.  
	\item Let $1 \leq i< n$ and $j=n+1$. Then $i$ and $j$ are separated by some $V_t$. 
	\item  $n$ and $n+1$ is separated by $V_{r+1}$. 
\end{enumerate} 
Thus $V'$ is a determining set for $K(n+1,k) $ and this proves the second inequality.\qed

\vspace{2 mm}

We next prove the following which, when put together with Theorem \ref{thm:rec_for_k_fixed}, will show that
when we fix $n$ and increase $k$, the determining number, as expected, either remains same or decrease.

{\theorem 
	\label{thm:rec_for_n_fixed}
	For positive integers $n,k$ with $n+1 \geq 2k+3$, $$Det(K(n+1,k+1))\leq Det(K(n,k)).$$}
\pf Let $Det(K(n,k))=r$ and $\{A_1,A_2,\ldots,A_r\}$ be a determining set for $K(n,k)$. Then $\cup_{i=1}^r A_i=[n]$ or $[n-1]$. However, by using techniques used in previous Theorem, without loss of generality, we can assume that $\cup_{i=1}^r A_i=[n-1]$. We set $V_i=A_i\cup \{n\}$ for $i=1,2,\ldots,r$. Then $|V_i|=k+1$. 

For $a,b \in \{1,2,\ldots,n-1\}$, there exists $A_i$ and hence $V_i$ which separates them. For $a \in \{1,2,\ldots,n-1\}$ and $b=n$, if $\cap_{i=1}^r A_i=\emptyset$, then there exists $A_i$ (and hence $V_i$) which separates $a$ and $b$. For $a \in \{1,2,\ldots,n-1,n\}$ and $b=n+1$, as $\cup_{i=1}^r V_i=[n]$, there exists $V_i$ which separates $a$ and $b$. 

So, let us assume that $\cap_{i=1}^r A_i\neq \emptyset$. However, as $\{A_1,A_2,\ldots,A_r\}$ be a determining set, $\cap_{i=1}^r A_i$ must be singleton, say $\{a\}$.

Thus, we have $\cap_{i=1}^r V_i=\{a,n\}$, $\cup_{i=1}^r V_i=[n]$ and $V_i$'s separate all pairs except $(a,n)$. We will modify $V_1$ to $V_1'$ such that $\{V_1',V_2,\ldots,V_r\}$ is a determining set for $K(n+1,k+1)$. We search for $t \in [n-1]\setminus V_1$ which can replace $a$ in $V_1$, as in that case, $V_1'$ will separate $a$ and $n$.

Suppose there does not exist any such $t$ such that $\{V_1'=V_1\cup\{t\}\setminus \{a\},V_2,V_3,\ldots,V_r\}$ is a determining set for $K(n+1,k+1)$. This implies that for every $t \in [n-1]\setminus V_1$, there exists $j_t\in V_1\setminus \{n\}$, i.e., $j_t \in A_1$ such that $(t,j_t)$ can not be separated by $\{V_1',V_2,\ldots,V_r\}$, i.e., $t,j_t$ are either both present or both absent in each of  $V_2,V_3,\ldots,V_r$ and $t,j_t \in V_1'$.

{\bf Claim:} If $t_1,t_2\in [n-1]\setminus V_1$ with $t_1\neq t_2$ such that they can not replace $a$ in $V_1$, then $j_{t_1}\neq j_{t_2}$.

{\bf Proof of Claim:} The proof goes in same line with that of proof of Claim 2 in previous theorem. However, for the sake of completeness, we write it down. Let there exist $t_1,t_2\in [n-1]\setminus V_1$ with $t_1\neq t_2$ such that they can not replace $a$ in $V_1$, but $j_{t_1}= j_{t_2}=j$ (say). Then $(t_1,j)$ and $(t_2,j)$ are either simultaneously present or absent in $V_2,V_3,\ldots,V_r$ and $t_1,t_2,j \in V_1'$. This implies that $t_1,t_2$ can not be separated by $V_1,V_2,\ldots,V_r$, a contradiction. Hence $j_{t_1}\neq j_{t_2}$ and the claim holds.

Now, there are $(n-1)-k$ elements in $[n-1]\setminus V_1$. If none of them can replace $a$ in $V_1$, then by the above Claim, there exists $n-k-1$ distinct elements in $V_1\setminus \{n\}$, i.e., $n-k-1\leq k$, i.e., $2k\geq n-1$. On the other hand, as $K(n+1,k+1)$ is a Kneser graph, we have $2k+3 \leq n+1$, that is, $2k \leq n-2 $. This is a contradiction.

Thus we can always find some $t \in [n-1]\setminus V_1$ which can replace $a$ in $V_1$ such that $\{V_1',V_2,\ldots,V_r\}$ is a determining set for $K(n+1,k+1)$. Thus, the theorem follows. \qed

\vspace{2 mm}

From Theorem \ref{thm:rec_for_k_fixed} and \ref{thm:rec_for_n_fixed}, we immediately have the following corollary. 

\begin{corollary}
	\label{cor:rec_for_n_fixed2}
	For positive integers $n,k$ with $n \geq 2k+3$, $$Det(K(n,k+1))\leq Det(K(n,k)).$$
\end{corollary}

\section{Determining Number of $K{(2k+1,k)}$}\label{exact-value-section}

\cite[Proposition 9]{boutin} used a linear algebraic construction to show the following result:
\begin{theorem}\label{thm:boutin_main_1}
	For any positive integer $r \geq 2$, we have 
	$Det(K(2^r-1,2^{r-1}-1))=r$.
\end{theorem}

These are precisely the points (shown in solid blue squares in Figure \ref{diagram}) and they lie on the straight line $n=2k+1$ and attain the lower bound $\lceil \log_2(n+1)\rceil$ described in Proposition  \ref{weak-bound}. In this section, we show that the lower bound is attained by all integer points on the line $n=2k+1$. 

Let $P[n,k]$ denote all the $k$-subsets of $[n]$, that is the vertices of $K(n,k)$. \begin{definition}
	\label{def:Auxiliary Set}
	For $n \geq 2k$, $\{V_1,V_2,\ldots,V_r\} \subset P[n,k]$ is called an auxiliary set if  $\displaystyle \cup_{i=1}^{r} V_i = [n]$ and 
	there exist no pair of distinct elements $a,b \in [n]$ such that  
	for each $i$, either 
	$\{a,b\} \subseteq V_i$ or $\{a,b\} \subseteq V_i^c$. 
	Define the auxiliary number of $K(n,k)$, denoted by $Aux(n,k)$, to be $r$ if $r$ is the 
	cardinality of a minimum auxiliary set. 
\end{definition}

Note that the definition of auxiliary set allows $n=2k$. From definition, it is clear that if $\{V_1,V_2,\ldots,V_r\} \subset P[n,k]$ 
is an auxiliary set for $K(n,k)$, then $\{V_1,V_2,\ldots,V_r\} \subset P[n,k]$ 
is a detemining set for $K(n,k)$ as well as $K(n+1,k)$.
Moreover, note that auxiliary sets require $\cup_{i=1}^r V_i=[n]$ which determining sets do not require. We start with the following Lemma which connects $Aux(2k,k)$ and $Aux(4k,2k)$.

{\lemma \label{lem:construction_of_Aux_Set_ofK2n2k}
	If $Aux(2k,k)=r$, then $Aux(4k,2k)\leq r+1$.} \\

\pf Let $\{A_1,A_2,\ldots,A_r\}$ be an auxiliary set for $K(2k,k) $. Let $A_i'=\{a+2k:a\in A_i\}$.
It is easy to see that $A_1',A_2',\ldots,A_r'$ separates any two elements in $[2k+1,4k]$ because if $A_i$ separates $(i,j)$ where $1 \leq i,j \leq 2k$, then $A_i'$ separates $(i+2k,j+2k)$ .  

For $i \in \{1,2,\ldots,r\}$, define $V_i=A_i \cup A_i'$ and $V_{r+1}=\{1,2,\ldots,2k\}$. We show that $\{V_1,V_2,\ldots,V_r,V_{r+1}\}$ is an auxiliary set for $K(4k,2k)$. 

Clearly, $\displaystyle \cup_{i=1}^{r+1} V_i=[4k]$. Consider the following cases.
\vspace{2 mm}

{ \bf Case 1: $1 \leq i,j \leq 2k$.}    In this case, $(i,j)$
is separated by some $A_t$ as $A_1,A_2,\ldots,A_r$ is an auxiliary 
set for $K(2k,k)$. Hence, $V_t$ separates the pair $(i,j)$. 

\vspace{ 2 mm}

{\bf Case 2: $2k+1 \leq i,j \leq 4k$.} The pair $(i+2k,j+2k)$
is separated by some $A_t'$ as mentioned earlier. Thus $V_t$ separates the pair $(i,j)$.  

\vspace{ 2 mm}

{\bf Case 3:} $1 \leq i \leq 2k$ and $2k+1 \leq j \leq 4k$ and $j-i \neq 2k$. Then $(i,j-2k)$ is separated by some $A_t$. Thus $(i,j)$ is separated by $V_t$. 
\vspace{ 2 mm}	

{\bf  Case 4:} $1\leq i \leq 2k$. In this case, $(i,i+2k)$
is separated by $V_{r+1}$. 

\vspace{ 2 mm}
Hence, all the pairs $(i,j)$ in $[4k]$ are separated by $\{V_1,V_2,\ldots,V_r,V_{r+1}\}$, proving the lemma. \qed

{\remark Note that, in the above proof, if $\cap_{i=1}^r A_i=\emptyset$, then $\cap_{i=1}^{r+1}V_i=\emptyset$.}

In the next lemma, we construct a determining set for $K(4k+3,2k+1)$, using a determining set for  $K(2k+1,k)$ with some particular properties. 

Before that, note if $\{V_1,V_2,\ldots,V_r\}$ is a set of vertices of $K(2k+1,k)$ such that $\cup _{i=1}^r V_i$ misses one point of $[2k+1]$, then without loss of generality, we can assume the missing point to be $2k+1$, i.e., $\cup _{i=1}^r V_i=[2k]$.

{\lemma	\label{lem:construction_of_Det_Set_ofK2n+12k}	Let $\{V_1,V_2,\ldots,V_r\}$ be a determining set for $K(2k+1,k)$ such that $\cap _{i=1}^r V_i = \emptyset$ and  $\cup _{i=1}^r V_i=[2k]$. Then, we can construct a determining set $\{W_1,W_2,\ldots,W_{r+1}\}$
	for $K(4k+3,2k+1)$ such that $\cap _{i=1}^{r+1} W_i = \emptyset$ and $\cup _{i=1}^{r+1} W_i=[4k+2]$. }\\
\\
\pf For $1 \leq i \leq r$, define  $V_i'=\{a+2k+1: a \in V_i\}$ and   $W_i= V_i \cup V_i' \cup \{4k+2\}$. Also define $W_{r+1} = \{1,2,\ldots,2k,2k+1\}$. 

We claim that $\{W_1,W_2,\ldots,W_{r+1}\}$ is a determining set for $K(4k+3,2k+1)$ with the above properties. 

Clearly, $\cup _{i=1}^r V_i'=[2k+2,4k+1]$ and hence $\cup _{i=1}^{r+1} W_i=[4k+2]$. Also $\cap _{i=1}^{r+1} W_i = \emptyset$. Now, consider the following cases.  
\begin{enumerate}
	\item Let $1 \leq i,j \leq 2k+1$. Then $i$ and $j$ are separated by some $V_t$ and hence by the corresponding $W_t$.
	\item Let $2k+2 \leq i,j \leq 4k+1$. Then $i$ and $j$ are separated by some $V_t'$ and hence by the corresponding $W_t$.
	\item Let $1 \leq i \leq 2k+1$ and $2k+2 \leq j \leq 4k+3$. Then $i$ and $j$ are separated $W_{r+1}$. 
	\item Let $2k+2 \leq i \leq 4k+1$  and $j=4k+2$. Since $\cap _{i=1}^r V_i = \emptyset$, we have $\cap _{i=1}^r V_i' = \emptyset$. Thus there exists $t \in \{1,2,\ldots,r\}$ such that $i \not\in V_t'$, i.e., $ i \not\in W_t$, but $4k+2 \in W_t$. Hence $i$ and $4k+2$ are separated $W_t$. 
	\item Let $1 \leq i \leq 4k+2$ and $j=4k+3$. Then $i$ and $j$ are separated by some $W_i$, for $i\in \{1,2,\ldots,r\}$.   
\end{enumerate} 
Combining all the cases, $\{W_1,W_2,\ldots,W_{r+1}\}$ is a determining set
for $K(4k+3,2k+1)$ with the aforesaid properties.\qed

{\proposition \label{prop:auxnumber}  
	For any positive integers $r$ and $k$ with $2^{r-1}-1 < 2k < 2^r-1$,  $Aux(2k,k)=r$ and there exists an auxiliary set $\{V_1,V_2,\ldots,V_r\}$ such that $\cap_{i=1}^r V_i=\emptyset$.} \\  
\\
\pf  We will prove this by induction on $r$. Our base case is $r=3$. For $r=3$, the permissible values of $k$ are $2$ and $3$, and we construct auxiliary sets of cardinality $3$ for each of $K(4,2)$ and $K(6,3)$.
\begin{enumerate} 
	\item $S=\{ V_1=\{1,2\}, V_2=\{1,3\}, V_3=\{2,4\} \}$ is an auxiliary set for $K(4,2)$. 	
	\item $T=\{ V_1=\{1,2,3\}, V_2=\{1,4,5\}, V_3=\{2,4,6\} \}$ is an auxiliary set for $K(6,3)$.	
\end{enumerate}
It can also be easily checked that $S$ and $T$ are auxiliary sets of minimum size for $K(4,2)$ and $K(6,3)$ respectively and $\cap_{i=1}^3 V_i=\emptyset$ in both the cases. Thus the result holds for $r=3$.

Now, we assume that for $r=t$ with all $k$ satisfying $2^{t-1}-1 < 2k < 2^t-1$, $Aux(2k,k)=t$ and there exists an auxiliary set $\{V_1,V_2,\ldots,V_t\}$ such that $\cap_{i=1}^t V_i=\emptyset$. 

Let $r=t+1$ and $k$ satisfy $2^{t}-1 < 2k < 2^{t+1}-1$. Consider the two following cases:

{\bf Case 1:} $k=2 \ell $ is even. Then we have $2^{t}-1 < 4 \ell  < 2^{t+1}-1$. So, $2 \ell < 2^t - \frac{1}{2} \implies 1/2 < 2^t-2 \ell \implies 1 < 2^t-2 \ell$. The last implication follows as $2$ divides the right hand side. So, finally we have $2^{t-1}-1 < 2 \ell < 2^{t}-1$. By induction hypothesis, $Aux(2 \ell, \ell)=t$ and there exists an auxiliary set $\{V_1,V_2,\ldots,V_t\}$ such that $\cap_{i=1}^t V_i=\emptyset$. Then by Lemma \ref{lem:construction_of_Aux_Set_ofK2n2k}  and the remark thereafter, $Aux(4 \ell,2 \ell)\leq t+1$ and there exists an auxiliary set $\{W_1,W_2,\ldots,W_{t+1}\}$ such that $\cap_{i=1}^{t+1} W_i=\emptyset$. On the other hand, by Proposition \ref{weak-bound}, $$Aux(4 \ell,2 \ell)  \geq \log_2(4 \ell+1)> \log_2(2^t)=t.$$ Thus $Aux(4 \ell,2 \ell)=Aux(2k,k)=t+1$ and there exists an auxiliary set $\{W_1,W_2,\ldots,W_{t+1}\}$ such that $\cap_{i=1}^{t+1} W_i=\emptyset$.

{\bf Case 2:} $k=2 \ell+1$ is odd.
Then we have $2^{t}-1 < 4\ell+2 < 2^{t+1}-1$, By similar arguments as in the previous case, we get $2^{t-1}-1 < 2 \ell < 2^t-1$. By induction hypothesis, $Aux(2 \ell, \ell)=t$ and there exists an auxiliary set $\{V_1,V_2,\ldots,V_t\}$ such that $\cap_{i=1}^t V_i=\emptyset$. Note that $\{V_1,V_2,\ldots,V_t\}$ is a determining set of $K(2 \ell+1, \ell)$ with $\cap_{i=1}^t V_i=\emptyset$ and $\cup_{i=1}^t V_i=[2 \ell]$. Then by Lemma \ref{lem:construction_of_Det_Set_ofK2n+12k}, there exists a determining set $\{W_1,W_2,\ldots,W_{t+1}\}$ of $K(4 \ell+3, 2\ell+1)$ such that $\cap_{i=1}^{t+1} W_i=\emptyset$ and $\cup_{i=1}^{t+1} W_i=[4 \ell+2]$. This implies that $\{W_1, W_2, \ldots, W_{t+1}\}$ is an auxiliary set for $K(4 \ell+2, 2\ell+1)$, i.e., $K(2k,k)$ such that $\cap_{i=1}^{t+1} W_i=\emptyset$. This means that $Aux(2k,k)\leq t+1$. Now by similar arguments as that in previous case, it can be shown that $Aux(2k,k)=t+1$.

Hence, by induction the proposition follows.\qed

We are now in a position to prove the main result of this section.

{\theorem \label{thm:main_theorem_sec2}$Det(K(2k+1,k))=r$ where $2^{r-1}-1 < 2k+1 \leq 2^r-1$, i.e., if $n=2k+1$, then $Det(K(n,k))=\left\lceil\log_2(n+1)\right\rceil$.}\\
\\
\pf By Proposition \ref{prop:auxnumber}, $Aux(2k,k)=r$, where $2^{r-1}-1 < 2k < 2^r-1$ and there exists an auxiliary set $\{V_1,V_2,\ldots,V_r\}$ such that $\cap_{i=1}^r V_i=\emptyset$. Thus $\{V_1,V_2,\ldots,V_r\}$ is a determining set for $K(2k+1,k)$ where $2^{r-1}-1 < 2k+1 \leq 2^r-1$, i.e., $Det(K(2k+1,k))\leq r$. 

Again, by Proposition 
\ref{weak-bound}, $Det(K(2k+1,k))\geq \log_2(2k+2)>\log_2(2^{r-1})=r-1$. Hence the theorem.  \qed

\begin{remark}
	One can see that Theorem \ref{thm:main_theorem_sec2} is an improvement of Theorem \ref{thm:boutin_main_1}, i.e., Theorem \ref{thm:boutin_main_1} shows that lower bound in Proposition \ref{weak-bound} is attained by some points on the line $n=2k+1$, whereas Theorem \ref{thm:main_theorem_sec2} establishes that the lower bound holds for all points on the line $n=2k+1$. 
\end{remark}

{\theorem \label{thm:n_equals_2k+2_andpwrof2} If $n=2k+2$ and $n$ is a power of $2$, then $Det(K(n,k))=\left\lceil\log_2(n+1)\right\rceil$. }\\
\\

\pf By Theorem \ref{thm:rec_for_k_fixed} and $Det(K(n,k))\geq \left\lceil \log_2(n+1)\right\rceil$, we have $$\left\lceil \log_2(2k+3)\right\rceil\leq Det(K(2k+2,k))\leq Det(K(2k+1,k))+1=\left\lceil\log_2(2k+2)\right\rceil +1.$$
Now, as $n=2k+2=2^s$, we have $k+1=2^{s-1}$ and $$\left\lceil \log_2(2k+3)\right\rceil=s+1=\left\lceil\log_2(2k+2)\right\rceil +1,$$
and hence the theorem follows.\qed 

{\remark The above theorem shows that we can find the exact value of determining number of Kneser graphs $K(n,k)$ for some integer points (shown in solid black squares in Figure \ref{diagram}) on the line $n=2k+2$. Note that these are precisely the corresponding points on the line $n=2k+1$ (shown in solid blue squares in Figure \ref{diagram}), for which exact values were determined by \cite{boutin}. }

{\corollary If $n=2k+2$, then $Det(K(n,k))=\left\lceil\log_2(n+1)\right\rceil$ or $\left\lceil\log_2(n+1)\right\rceil+1$.} 

\begin{proof}
	If $n=2k+2$ is a power of $2$, by Theorem 
	\ref{thm:n_equals_2k+2_andpwrof2}, we have the result. When $n=2k+2$ is not a power of $2$, 
	using Theorem \ref{thm:rec_for_k_fixed} we get   $$ Det(K(2k+2,k))\leq Det(K(2k+1,k))+1=\left\lceil\log_2(2k+2)\right\rceil +1\leq \left\lceil\log_2(n+1)\right\rceil +1.$$
	Again from Proposition \ref{weak-bound}, we have $ Det(K(2k+2,k))\geq \left\lceil\log_2(n+1)\right\rceil$. Thus when $n=2k+2$ is not a power of $2$, we proved that $Det(K(2k+2,k))$ lies between two consecutive integers and hence has to be one of them, thereby completing the proof. 
\end{proof}

\section{Bounds for $Det(K(n,k))$}\label{Bounds-section}

\subsection{Lower Bound}
\cite{boutin} proved a lower bound for the determining number of any Kneser graph $K(n,k)$ which was mentioned in Proposition \ref{weak-bound}. Next, we provide a lower bound which is stronger than that in \cite{boutin} if $k=\Omega(n/\log n)$. 

{\theorem
	\label{thn:lowerbound}
	For any positive integers $n,k$ with $k < \frac{n}{2}$, 
	$Det(K(n,k)) \geq \frac{2n-2}{k+1} $. } \\

\pf
Let $Det(K(n,k))=r$ and let $\{V_1,V_2,\ldots,V_r\}$ be a determining set for $K(n,k) $. Hence, $\sum_{i=1}^r|V_i| = rk$. Now, we count $\sum_{i=1}^r|V_i|$  in another way. First, we note the following:

\begin{enumerate}
	\item $\displaystyle \cup_{i=1}^{r} V_i$ can miss at most one element of $[n]$.
	
	\item There can be at most $r$ elements which occur in exactly one of the sets. If not, suppose there are $t$ ($t >r$) elements which occur in exactly one set.  In this case at least two of these $t$ elements (say $a$ and $b$) would occur in the same set and they occur only in that set. Hence, $a$ and $b$ are not separated which is a contradiction. 
\end{enumerate}
Thus, all other elements of $[n]$ are there in at least two sets.  Hence, counting according to the number of appearances of any element in $[n]$, yields the following equation 

$$
r+2(n-1-r)  \leq   rk 
\implies 2n-2-r  \leq rk  
\implies r \geq \frac{2n-2}{k+1},
$$
completing the proof.  \qed

\subsection{Upper Bound}

We are now interested in constructing an improved upper bound for $D(K(n,k))$. \cite{kneser} proved the following theorem.

{\theorem (\cite[Theorem 3.1]{kneser})
	\label{thm:Caceres}
	For positive integers $n,k$ with $2k \leq n \leq  \frac{k(k+1)}{2}$, 
	$Det(K(n,k)) \leq k$. 
}

Here, we provide a stronger upper bound by using an explicit construction of a determining set. 

{\theorem
	\label{thm:stronger_upper_bound}
	Let $n,k,r$ be positive integers with $n \leq k(k+1)/2$
	and $k\geq  r\geq 3$. Then, for all integers $n$ with  $ n \leq r(r+1)/2+1 \implies Det(K(n,k)) \leq  r$. 
} 

\pf At first we prove that  $Det(K(n,r)) \leq  r$ when $n= r(r+1)/2+1$.  Consider the following $r$-sets
\begin{eqnarray}
\label{eqn:constructionofVi}
V_1 & := &\{1,2,3,\ldots,r-1,r\} , \nonumber  \\
V_2 & := &\{1,r+1,r+2,\ldots,2r-2,2r-1\},  \nonumber \\
V_3 & := &\{2,r+1,2r,\ldots,3r-4,3r-3\},  \nonumber \\
V_4 & := &\{3,r+2,2r,\ldots,4r-7,4r-6\},  \nonumber \\  
\hdots & & \hdots \nonumber \\
V_r & :=& \{r-1,2r-2,3r-4,\ldots,(r-1)(r+2)/2,r(r+1)/2\}.  \nonumber 
\end{eqnarray}

Our main idea behind the construction is as follows: After $V_1,V_2, \dots, V_{i-1}$ being already constructed, we construct $V_i$ as follows:
\begin{multline*}
V_i:= \displaystyle \{V_1^{i-1},\dots,V_{i-1}^{i-1},r(i-1)-\frac{(i-1)(i-2)}{2}+1, \ r(i-1)-\frac{(i-1)(i-2)}{2}+2, \dots , \\
r(i-1)-\frac{(i-1)(i-2)}{2}+ r-i+1 \}.
\end{multline*}

Here, $V_i^j$ denotes the $j$-th element of $V_i$ when we write all the elements of $V_i$ in ascending order.
We claim the following:\\
{\bf Claim:} $V=\{V_1,V_2,\ldots,V_r\}$  is a determining set for $K(n,k)$. \\
{\bf Proof of Claim:}
Let $S=\{r,2r-1,3r-3,\ldots,\frac{r(r+1)}{2}\}$ and $V_i'=V_i \setminus S$. We observe the following. 

\begin{enumerate} 
	
	\item Each element of $S$ is exactly in one of the $V_i$ and each element of $[n]\setminus S$ is in exactly two of the $V_i$'s. Hence any $a \in S$ is separated from any $b \in [n]\setminus S$. Also any two elements in $S$ are separated by some $V_i$.
	
	\item Already having constructed $V_1, V_2, \ldots, V_{i-1}$, we construct $V_i$ in a way such that $|V_i \cap V_j| = 1$ for all $ j < i$. Hence,  by our construction itself, $|V_i \cap V_j|=1$ for all $i,j$. Besides, each element $a \in [n]\setminus S$ is in exactly two sets. Hence any two elements $a,b \in [n]\setminus S$ are separated. 
	
\end{enumerate}
Thus,  $V=\{V_1,V_2,\ldots,V_r\}$  is a determining set for $K(n,k)$ and the claim holds. Hence, $Det(K(n,r)) \leq  r$ when $n= r(r+1)/2+1$. 

Now, Theorem \ref{thm:rec_for_n_fixed} states that the sequence $Det(K(n,r))$ is weakly decreasing when we keep $n$ fixed and increase $r$. Thus, for
$n=r(r+1)/2 + 1$, we have $Det(K(n,k)) \leq Det(K(n,r)) \leq r$. Now for $n$ with $n \leq r(r+1)/2+1 $,  we use Theorem \ref{thm:rec_for_k_fixed} directly to get $Det(K(n,k)) \leq r$. \qed

\begin{remark}
	When $r\ll k$,  Theorem \ref{thm:stronger_upper_bound} clearly gives a much better bound than Theorem \ref{thm:Caceres}. When $r=k$,  this is Theorem \ref{thm:Caceres}. In Figure \ref{diagram}, the yellow region is where the exact value of determining number is known from Theorem \ref{caceres_first_main_result}. Caceres {\it et.al.} in Theorem \ref{thm:Caceres} also showed that above the yellow region, $Det(K(n,k))\leq k$. Using Theorem \ref{thm:stronger_upper_bound}, we prove a stronger upper bound as shown in different shades of gray in Figure \ref{diagram}.
\end{remark}

\acknowledgements
\label{sec:ack}
The authors are grateful to the anonymous referees for the useful comments. The authors would also like to thank Niranjan Balachandran for suggesting the problem and Arnab Mandal for several fruitful discussions in the initial phase of this work.

\nocite{*}
\bibliographystyle{abbrvnat}
\bibliography{sample-dmtcs}

\begin{thebibliography}{9}
\providecommand{\natexlab}[1]{#1}
\providecommand{\url}[1]{\texttt{#1}}
\expandafter\ifx\csname urlstyle\endcsname\relax
  \providecommand{\doi}[1]{doi: #1}\else
  \providecommand{\doi}{doi: \begingroup \urlstyle{rm}\Url}\fi

\bibitem[Albertson and Boutin(2007)]{boutin-distinguish}
M.~O. Albertson and D.~L. Boutin.
\newblock Using determining sets to distinguish {K}neser graphs.
\newblock \emph{The Electronic Journal of Combinatorics}, 14\penalty0
  (1):\penalty0 R20, 2007.

\bibitem[Boutin(2006)]{boutin}
D.~L. Boutin.
\newblock Identifying graph automorphisms using determining sets.
\newblock \emph{The Electronic Journal of Combinatorics}, 13:\penalty0 R78,
  2006.

\bibitem[Boutin(2009)]{boutin-product}
D.~L. Boutin.
\newblock The determining number of a {C}artesian product.
\newblock \emph{Journal of Graph Theory}, 61\penalty0 (2):\penalty0 77--87,
  2009.

\bibitem[C{\'a}ceres et~al.(2013)C{\'a}ceres, Garijo, Gonz{\'a}lez,
  M{\'a}rquez, and Puertas]{kneser}
J.~C{\'a}ceres, D.~Garijo, A.~Gonz{\'a}lez, A.~M{\'a}rquez, and M.~L. Puertas.
\newblock The determining number of {K}neser graphs.
\newblock \emph{Discrete Mathematics and Theoretical Computer Science},
  15\penalty0 (1):\penalty0 1--14, 2013.

\bibitem[Das(2020)]{das-caldam}
A.~Das.
\newblock Determining number of generalized and double generalized {P}etersen
  graph.
\newblock In \emph{Conference on Algorithms and Discrete Applied Mathematics},
  pages 131--140. Springer, 2020.

\bibitem[Das and Saha(2020)]{das-saha}
A.~Das and M.~Saha.
\newblock Determining number of some families of cubic graphs.
\newblock \emph{Journal of Algebra and Related Topics}, 8\penalty0
  (2):\penalty0 39--55, 2020.

\bibitem[Erwin and Harary(2006)]{harary}
D.~Erwin and F.~Harary.
\newblock Destroying automorphisms by fixing nodes.
\newblock \emph{Discrete Mathematics}, 306\penalty0 (24):\penalty0 3244--3252,
  2006.

\bibitem[Godsil and Royle(2001)]{godsil-royle}
C.~Godsil and G.~F. Royle.
\newblock \emph{Algebraic graph theory}, volume 207.
\newblock Springer Science \& Business Media, 2001.

\bibitem[Pan and Guo(2019)]{coprime}
J.~Pan and X.~Guo.
\newblock The full automorphism groups, determining sets and resolving sets of
  coprime graphs.
\newblock \emph{Graphs and Combinatorics}, 35\penalty0 (2):\penalty0 485--501,
  2019.

\end{thebibliography}
\label{sec:biblio}

\end{document}